\documentclass[12pt, twoside, a4paper, fleqn]{article}

\usepackage{latexsym}
\usepackage{amscd}
\usepackage{theorem}
\usepackage{pifont}
\usepackage{mathbbol}
\usepackage{amsfonts}
\usepackage{xspace}
\usepackage{amssymb}
\usepackage{fancyhdr}
\usepackage{mathrsfs}
\usepackage{amsmath}%
\usepackage{graphics}%
\usepackage{graphicx}%
\usepackage{euscript}%
\usepackage{psfrag}%
\usepackage{empheq}%
\usepackage{upgreek}
\usepackage{eufrak}

\DeclareMathAlphabet{\mathpzc}{OT1}{pzc}{m}{it}

{\theorembodyfont{\slshape} \newtheorem{theorem}{Theorem}[section]}
{\theorembodyfont{\slshape} }
{\theorembodyfont{\slshape} }
{\theorembodyfont{\slshape} \newtheorem{corollary}[theorem]{Corollary}}

\newenvironment{proof}{\noindent\textbf{Proof:\ }}{$\hfill{\bullet}$}

\numberwithin{equation}{section}

\bigskip

\title{Lyapunov exponents of polynomials with respect to certain weighted Lyubich's measures}
\bigskip
\bigskip

\author{{\sc Shrihari Sridharan\footnote{This author was supported by a Fasttrack Grant for Young Scientists awarded by the Department of Science and Technology, Government of India, vide SR/FTP/MS-008/2012.},\ \ \ \ \ \ Atma Ram Tiwari} \\ Indian Institute of Science Education and Research \\ Thiruvananthapuram (IISER-TVM) \\ {\tt shrihari@iisertvm.ac.in},\ \ \ \ \ \ {\tt artiwari15@iisertvm.ac.in}}
\bigskip

\date{}
\bigskip 

\begin{document}

\maketitle

\begin{abstract}
\noindent 
In this paper, we consider a monic, centred, hyperbolic polynomial of degree $d \ge 2$, restricted on its Julia set and compute its Lyapunov exponents with respect to certain weighted Lyubich's measures. In particular, we show a certain well-behavedness of some coefficients of the Lyapunov exponents, that quantifies the non-well-behavedness in a system. 
\end{abstract}

\begin{tabular}{l c l}
\textbf{Keywords} & : & Lyapunov exponents, \\
& & Weighted Lyubich's measures, \\ 
& & Pressure function \\ 
\\ 
\textbf{AMS Subject Classifications} & : & 37B25, 37F15, 37F10. 
\end{tabular}
\bigskip  

\thispagestyle{empty}


\section{Introduction}

\noindent 
The theory of dynamical systems is studied widely by scientists from all disciplines, owing to their applications in chaos theory, bifurcation theory, scientific modelling, image processing \textit{etc}. A major characteristic computable quantity that describes the well-behavedness or otherwise of dynamical systems is the Lyapunov exponent. The Lyapunov exponent indicates the rate of separation of typical trajectories that begin at ponits which are sufficiently close to each other, thus, quantitatively describing the senstivity of the chaotic dynamical system to the parameters involved. 
\medskip 

\noindent 
The relationship between the Lyapunov exponent and exponential divergence of typical trajectories, as studied by Benettin \textit{et al} in \cite{Benettin}, have been widely used to study various dynamical systems. This study has been continued by various people in different fields of mathematics as well as physics. The relationship between Lyapunov exponent and other computable quantities such as topological entropy, index of exponential stability and supersymmetry are well-explored topics amidst mathematicians and physicists, for example, \cite{Benettin, tc:19, bcnv:18, gr:94, Pesin, Ruelle, Sridharan, Zinsmeister}. In particular, the relationship between the Lyapunov exponent and the degree of stochasticity for typical trajectories like Kolmogorov entropy has been calculated by Pesin in \cite{Pesin}. 
\medskip 

\noindent 
The concept of using symbolic spaces is a powerful tool in calculating the Lyapunov exponent for any dynamical system. This has been used in the study of dynamical systems, by various authors, including \cite{Coelho, Lyubich, Pesin, Ruelle, Sridharan, Zinsmeister}. The authors, in \cite{SridharanArt}, make use of this concept and compute the Lyapunov exponent of quadratic and cubic polynomials satisfying a few assumptions, with respect to various measures associated to the family of Bernoulli measures on appropriate symbolic spaces. In fact, they establish the dependence of the Lyapunov exponent on the coefficients of the polynomials and their relationship to the Hausdorff dimension of the respective Julia set. 
\medskip 

\noindent 
A surprising, yet powerful result by the authors in \cite{SridharanArt} was regarding the second derivative of the Lyapunov exponent of a monic, centred, hyperbolic quadratic polynomial parametrised by a complex co-efficient restricted on its Julia set, with respect to the real part and the imaginary part of the parameter, as well as the second derivative of the Lyapunov exponent of a monic, centred, hyperbolic cubic polynomial parametrised by two complex coefficients restricted on its Julia set, with respect to the respective real part and the imaginary part of the same parameter agreed in size, but varied in sign. This was captured as the second point in the final section on concluding observations. This paper naturally grows out of that observation to verify the same for a polynomial with similar properties, however of any arbitrary degree, say $d$. We will now state the main theorems of this paper.
\medskip 

\noindent
\begin{theorem} 
\label{theorem1}
Consider a monic, centred, hyperbolic polynomial $P$, of degree $d > 1$, with complex coefficients, given by 
\[ P(z)\ \ :=\ \ z^{d} + \left( \alpha_{d - 2} + i \beta_{d - 2} \right) z^{d - 2} + \cdots + \left( \alpha_{1} + i \beta_{1} \right) z + \left( \alpha_{0} + i \beta_{0} \right), \] 
where $\alpha_{r} \in \mathbb{R}$ and $\beta_{r} \in \mathbb{R}$ satisfying $\alpha_{r}^{2} + \beta_{r}^{2} < 1$ for all $0 \le r \le d - 2$, restricted on its Julia set, $\mathcal{J}_{P}$ with bounded critical orbit. Suppose $\Lambda_{\mu}$ is the Lyapunov exponent of $P$ with respect to a probability measure $\mu$ supported on $\mathcal{J}_{P}$ given by 
\begin{equation} 
\label{Lyapdef}
 \Lambda_{\mu} (P)\ \ :=\ \ - \int_{\mathcal{J}_{P}} \log |P'| d \mu.
 \end{equation} 
Then, irrespective of the strictly positive probability vector $\vec{p} = (p_{1}, p_{2}, \cdots, p_{d})$ based on which we define a probability measure $\mu_{\vec{p}}$ supported on $\mathcal{J}_{P}$, we have 
\[ \frac{\partial^{2} \Lambda_{\mu_{\vec{p}}}}{\partial \alpha_{r} \partial \alpha_{s}}\ \ =\ \ - \frac{\partial^{2} \Lambda_{\mu_{\vec{p}}}}{\partial \beta_{r} \partial \beta_{s}},\ \ \ \ \text{for every}\ \ 0 \le r, s \le d - 2. \] 
\end{theorem} 
\medskip 

\noindent 
In the course of proving this theorem, we will evaluate the Lyapunov exponent of a monic, centred, hyperbolic polynomial $P$ of degree $d$, up to terms with a certain order. Then, we also obtain the following theorem pertaining to the first observation in \cite{SridharanArt}. 
\medskip 

\noindent 
\begin{theorem} 
\label{theorem2}
Consider the polynomial $P$ with bounded critical orbit; satisfying the hypothesis in theorem \eqref{theorem1} with real coefficients and complex coefficients separately as given below: 
\begin{eqnarray} 
\label{polyreal} 
P_{\mathbb{R}} (z) & = & z^{d} + \alpha_{d - 2} z^{d - 2} + \cdots + \alpha_{1} z + \alpha_{0}, \\ 
\label{polycomp} 
P_{\mathbb{C}} (z) & = & z^{d} + \left( \alpha_{d - 2} + i \beta_{d - 2} \right) z^{d - 2} + \cdots + \left( \alpha_{1} + i \beta_{1} \right) z + \left( \alpha_{0} + i \beta_{0} \right). 
\end{eqnarray} 
Let $\Lambda_{\mu_{\vec{p}}}$ denote the appropriate Lyapunov exponent with respect to the probability measure $\mu_{\vec{p}}$ based on the strictly positive probability vector $\vec{p} = (p_{1}, p_{2}, \cdots, p_{d})$. Then, 
\[ \frac{\partial \Lambda_{\mu_{\vec{p}}} (P_{\mathbb{R}})}{\partial \alpha_{r}}\ \ =\ \ \frac{\partial \Lambda_{\mu_{\vec{p}}} (P_{\mathbb{C}})}{\partial \alpha_{r}}; \qquad \qquad \qquad \frac{\partial^{2} \Lambda_{\mu_{\vec{p}}} (P_{\mathbb{R}})}{\partial \alpha_{r}^{2}}\ \ =\ \ \frac{\partial^{2} \Lambda_{\mu_{\vec{p}}} (P_{\mathbb{C}})}{\partial \alpha_{r}^{2}}. \]
\end{theorem}

\noindent
This paper is structured as follows: In the following section namely section \eqref{setting}, we state the necessary definitions and build the necessary preliminaries that would give us the basic settings on which the results of this paper rest. In section \eqref{weight}, we define the equidistributed Lyubich's measure and the weighted Lyubich's measure, by according different weights on the preimage branches of generic points under the considered polynomial map. In section \eqref{cofn}, we make use of a topological conjugacy and calculate certain coefficient functions (only as much as necessary), that will come in handy while we compute the Lyapunov exponents in section \eqref{lecom}. With all the computations done, we prove our main theorems in section \eqref{proof} and conclude the paper, with corollaries by interpreting the theorem for the specific cases when the polynomial is quadratic or cubic. 

\section{Basic settings} 
\label{setting}

\noindent 
In this section, we focus on building the basic settings of this paper by writing the essential terminologies, elementary definitions and a few results from the literature that we will require to prove our theorems. 
\medskip 

\noindent 
We consider a monic centred polynomial, say $P$ of degree $d > 1$, with complex coefficients defined on the Riemann sphere $\overline{\mathbb{C}} = \mathbb{C} \cup \{ \infty \}$, given by 
\[ P(z)\ \ =\ \ A_{d} z^{d} + A_{d - 1} z^{d - 1} + A_{d - 2} z^{d - 2} + \cdots + A_{1} z + A_{0}. \]
Here, monic means $A_{d} \equiv 1$ and centred means $A_{d - 1} \equiv 0$. Further, we consider the remaining coefficients to be complex; $A_{r} = \alpha_{r} + i \beta_{r}$ with $\alpha_{r}, \beta_{r} \in \mathbb{R}$ for $0 \le r \le d - 2$. It is then obvious that the critical point of $P$ is determined by the $(d - 2)$ coefficients $A_{d - 2}, \cdots, A_{1}$ while its critical orbit is determined by the $(d - 1)$ coefficients $A_{d - 2}, \cdots, A_{1}, A_{0}$. 
\medskip 

\noindent 
In this paper, we shall be interested in restricting the parameter space of the polynomial $(A_{d - 2}, \cdots, A_{1}, A_{0}) \in \mathbb{C}^{d - 1}$ so that the critical orbit of $P$ remains bounded. For example, the origin in $\mathbb{C}^{d - 1}$ is a point that yields the polynomial $P(z) = z^{d}$, whose critical point is the fixed point of the polynomial at $0 \in \mathbb{C}$, and thus has a bounded critical orbit. We reserve the notation $Q(z)$ to represent this polynomial, throughout this paper, where all, but the leading coefficient, are $0$; \textit{i.e.}, $Q(z) \equiv z^{d}$. 
\medskip 

\noindent 
For such a polynomial $P$ with bounded critical orbits, we define its filled Julia set $\mathcal{K}_{P}$ and its basin of attraction to the point at infinity $\mathcal{A}_{P} (\infty)$ respectively as 
\[ \mathcal{K}_{P}\ \ :=\ \ \left\{ z \in \mathbb{C} : P^{n} (z) \nrightarrow \infty\ \text{for any}\ n \right\}; \qquad \qquad \mathcal{A}_{P} (\infty)\ \ :=\ \ \left\{ z \in \mathbb{C} : P^{n} (z) \rightarrow \infty \right\}. \] 
It is, of course, obvious from the definition that the completely $P$-invariant sets $\mathcal{K}_{P}$ and $\mathcal{A}_{P} (\infty)$ dichotomise the Riemann sphere, $\overline{\mathbb{C}}$. The common topological boundary between the sets $\mathcal{K}_{P}$ and $\mathcal{A}_{P} (\infty)$ is then defined to be the \emph{Julia set} of the polynomial, denoted by $\mathcal{J}_{P}$. Interested readers may know that the Julia set of a polynomial map $P$ of degree $d > 1$ is defined in various ways. $\mathcal{J}_{P}$ is the closure of the set of all periodic points that satisfy a repelling condition, 
\[ \mathcal{J}_{P}\ \ =\ \ \overline{\left\{ z_{0} \in \mathbb{C} : P^{m} z_{0} = z_{0}\ \text{for some}\ m \in \mathbb{Z}_{+}\ \text{and}\ \left| \left( P^{m} \right)' (z_{0}) \right| > 1 \right\}}. \]
Alternatively, $\mathcal{J}_{P}$ is the set of points where the family of iterates of $P$, \textit{i.e.}, $\left\{ P^{n} \right\}_{n \ge 1}$ does not form a normal family (in the sense of Montel). The various definitions elucidate that $\mathcal{J}_{P}$ is a non-empty, compact, completely $P$-invariant metric space. Observe that the definitions entail, for $Q(z) = z^{d}$, we have 
\[ \mathcal{K}_{Q}\ =\ \left\{ z \in \mathbb{C} : |z| \le 1 \right\}; \qquad \mathcal{A}_{Q} (\infty)\ =\ \left\{ z \in \mathbb{C} : |z| > 1 \right\}; \qquad \mathcal{J}_{Q}\ =\ \left\{ z \in \mathbb{C} : |z| = 1 \right\}. \]
For more properties of the Julia set, one may refer \cite{Beardon, Lyubich}. 
\medskip 

\noindent 
The rationale behind considering only the set of all monic, centred polynomial maps $P$, given by 
\begin{equation} 
\label{polynomial}
P(z)\ \ =\ \ z^{d} + A_{d - 2} z^{d - 2} + \cdots + A_{1} z + A_{0}, 
\end{equation} 
with coefficients $A_{r} = \alpha_{r} + i \beta_{r}$ is due to the fact that any polynomial map $\mathcal{P}$ of degree $d$ given by 
\[ \mathcal{P} (z)\ \ =\ \ B_{d} z^{d} + B_{d - 1} z^{d - 1} + \cdots + B_{1} z + B_{0}, \] 
can, by an affine change of coordinates, be written as in equation \eqref{polynomial}. For the sake of computations in this paper, we demand that the coefficients $A_{r}$ in equation \eqref{polynomial} satisfy $|A_{r}|^{2} = \alpha_{r}^{2} + \beta_{r}^{2} < 1$ for $d - 2 \ge r \ge 0$. We concentrate only on hyperbolic polynomials $P$ restricted on their respective Julia set $\mathcal{J}_{P}$. Here, by \emph{hyperbolicity}, we mean that there exists constants $C > 0$ and $\lambda > 1$ such that for any $z \in \mathcal{J}_{P}$, we have $| P^{n} (z) | \ge C \lambda^{n}$ for all $n \ge 1$. It is an easy observation that a hyperbolic $\mathcal{J}_{P}$ is topologically connected.  
\medskip 

\noindent 
By a result of Lyubich, as in \cite{Lyubich}, the family of hyperbolic polynomial maps of the same degree is structurally stable. This means that when the hyperbolic polynomials $P_{1}$ and $P_{2}$ are restricted on their respective Julia sets $\mathcal{J}_{P_{1}}$ and $\mathcal{J}_{P_{2}}$, there exists a conjugacy, say $\Phi_{\left( P_{1}, P_{2} \right)} : \mathcal{J}_{P_{1}} \longrightarrow \mathcal{J}_{P_{2}}$ that satisfies $\Phi_{\left( P_{1}, P_{2} \right)} \circ P_{1} = P_{2} \circ \Phi_{\left( P_{1}, P_{2} \right)}$. The conjugacy $\Phi_{\left( P_{1}, P_{2} \right)}$ is naturally dependent on the polynomials $P_{1}$ and $P_{2}$, in other words the respective coefficients. In particular, we shall be interested in the conjugacy between the unit circle $\mathbb{S}^{1}$ which is the Julia set of $Q(z) = z^{d}$ and the Julia set of the monic centred hyperbolic polynomial map $P$, namely $\Phi_{\left( Q, P \right)} : \mathbb{S}^{1} = \mathcal{J}_{Q} \longrightarrow \mathcal{J}_{P}$ as written in equation \eqref{polynomial}, namely,
\begin{equation} 
\label{conjugacy}
\Phi_{\left( Q, P \right)} \left( z^{d} \right) - P \left( \Phi_{\left( Q, P \right)} (z) \right)\ \ =\ \ 0. 
\end{equation} 
Here, we keep the polynomial $Q$ fixed and focus on the dependence of the conjugacy on the coefficients of $P$. Thus, we suppress $Q$ in the notation of the conjugacy and denote the same merely by $\Phi_{P}$. It is then a result from \cite{Zinsmeister, Lyubich, cg:93} that has been interpreted in theorems (4.1) and (4.2) in \cite{SridharanArt} that the conjugacy $\Phi_{P}$ depends analytically on each of the parameters, $A_{r}$ for $d - 2 \ge r \ge 0$. 

\section{Weighted Lyubich's measures} 
\label{weight}

\noindent 
In this section, we shall focus on the space of $P$-invariant probability measures supported on the compact metric space $\mathcal{J}_{P}$, denoted by $\mathcal{M}_{P}$, \textit{i.e.}, 
\[ \mathcal{M}_{P}\ \ :=\ \ \left\{ \mu : \mathcal{J}_{P} \longrightarrow [0, 1]\ :\ \mu(E) = \mu(P^{-1} E)\ \forall E \subseteq \mathcal{J}_{P} \right\}. \] 

\noindent 
Owing to the density of preimages of any generic point $\zeta \in \mathcal{J}_{P}$, it can be observed that the sequence of measures 
\begin{equation} 
\label{lyubich}
\mu_{n}^{(\zeta)}\ \ :=\ \ \frac{1}{d^{n}} \sum_{P^{n} \omega\, =\, \zeta} \delta_{\omega},\ \ \text{where}\ \delta_{\omega}\ \text{is the Dirac delta measure at the point}\ \omega, 
\end{equation} 
converges to some measure $\mu \in \mathcal{M}_{P}$ called the \emph{Lyubich's measure}, independent of $\zeta$, in the weak*-topology, see for example \cite{Steinmetz}. For example, the polynomial map $P(z) = z^{d}$ has the unit circle, $\mathbb{S}^{1}$ in $\mathbb{C}$ as its Julia set and the Lyubich's measure can be thought of as the Haar measure on $\mathbb{S}^{1}$. It is then obvious that the support of the Lyubich's measure is the Julia set. 
\medskip 

\noindent 
An effective implication of the Lyubich's measure is that if the hyperbolic $\mathcal{J}_{P}$ is divided into $d$ mutually disjoint equal arcs, then every arc contains one and only one preimage of any generic point $\zeta \in \mathcal{J}_{P}$. Thus, by according equal weightage to every preimage branch, we obtain an equilibrium distribution. However, suppose we accord different weights to the different preimage branches of $P$, we obtain a distorted distribution, as we now explain. 
\medskip 

\noindent 
Let $\left\{ P_{1}, P_{2}, \cdots, P_{d} \right\}$ denote the preimage branches of the polynomial map $P$. For some strictly positive probability vector, $\vec{p}$, \textit{i.e.}, $\vec{p} = (p_{1}, p_{2}, \cdots, p_{d})$ with $p_{j} > 0,\ \forall j$ and $\sum\limits_{j = 1}^{d} p_{j} = 1$, we define a quantity called the weighted Lyubich's measure as follows. Consider any $n$-lettered word $\eta = \left( \eta_{1}, \eta_{2}, \cdots, \eta_{n} \right) \in \left\{ 1, 2, \cdots, d \right\}^{n}$. We define the appropriate $n$-th order preimage branch of $P$ as $P_{\eta} = P_{\eta_{n}} \circ P_{\eta_{n - 1}} \circ \cdots \circ P_{\eta_{1}}$. Then, for any generic point $\zeta \in \mathcal{J}_{P}$, we define a sequence of measures with respect to $\vec{p}$ namely $\left\{ \left( \mu_{\vec{p}}^{(\zeta)} \right)_{n} \right\}_{n\, \ge\, 1}$ as 
\begin{equation} 
\label{weightedlyubich} 
\left( \mu_{\vec{p}}^{(\zeta)} \right)_{n}\ \ :=\ \ \sum_{\eta\ :\ P_{\eta} \omega\, =\, \zeta} p_{\eta_{n}} p_{\eta_{n - 1}} \cdots p_{\eta_{1}} \delta_{\omega}. 
\end{equation} 
Then, owing to the uniform distribution of the preimage branches in $\mathcal{J}_{P}$, when $p_{j}$'s are not equally distributed in the probability vector $\vec{p}$, some preimage branches gain prominence over the other ones. However, as $n$ increases, the sections of the Julia set, $\mathcal{J}_{P}$, that gain prominence get tinier and tinier. In any case, this sequence of measures $\left\{ \left( \mu_{\vec{p}}^{(\zeta)} \right)_{n} \right\}_{n\, \ge\, 1}$ converges to some measure $\mu_{\vec{p}} \in \mathcal{M}_{P}$ called the \emph{weighted Lyubich's measure}, independent of $\zeta$, in the weak*-topology. 
\medskip 

\noindent 
This family of weighted Lyubich's measures is interesting to work with, especially when one of the $p_{j}$'s is extremely close to $1$, leaving the remainder of the $p_{k}$'s to be arbitrarily close to $0$. In such a case, the section of the Julia set corresponding to the preimage branch that gains prominence in the Julia set, eventually reduces to a point measure. 
\medskip 

\noindent 
We conclude this section by defining a quantity called pressure for a real-valued continuous function, say $f$, defined on $\mathcal{J}_{P}$, in accordance with thermodynamic formalism. This will come in handy, when we compute the Lyapunov exponent, later. 
\begin{equation}
\label{pressure} 
\mathfrak{P} (f)\ \ :=\ \ \sup_{\mu\, \in\, \mathcal{M}_{P}} \left\{ h_{\mu} (P) + \int_{\mathcal{J}_{P}} f d \mu \right\}, 
\end{equation} 
where $h_{\mu} (P)$ represents the measure theoretic entropy of $P$ with respect to the measure $\mu$. For more properties of pressure and entropy, interested readers are referred to \cite{Walter}. 
\medskip 

\noindent 
In this paper, we are interested in the real-valued continuous function $- \log |P'|$ defined on $\mathcal{J}_{P}$. Hence, 
\begin{eqnarray} 
\label{pratwork}
\mathfrak{P} \left( - \log |P'| \right) & = & \sup_{\mu\, \in\, \mathcal{M}_{P}} \left\{ h_{\mu} (P) - \int_{\mathcal{J}_{P}} \log |P'| d \mu \right\} \nonumber \\ 
& = & \sup_{\mu\, \in\, \mathcal{M}_{P}} \left\{ h_{\mu} (P) + \Lambda_{\mu} (P) \right\} \nonumber \\
& = & \sup_{\nu\, \in\, \mathcal{M}_{Q}} \left\{ h_{\nu} (Q) - \log d - \int_{\mathbb{S}^{1}} \log | \Phi (z) | d \nu \right\}. 
\end{eqnarray} 

\section{Coefficient functions} 
\label{cofn}

In this section, we make necessary preparations to calculate the Lyapunov exponent, as defined in equation \eqref{Lyapdef} of a monic, centred, hyperbolic polynomial $P$, as given in equation \eqref{polynomial}, that will help us obtain expressions, as required in the main theorems. Since we know from section \eqref{setting} that the conjugacy $\Phi_{P}$ analytically depends on the coefficients $A_{r}$ for $d - 2 \ge r \ge 0$, it is only reasonable to consider 
\[ \Phi_{P} (z)\ \ =\ \ z + \sum_{\xi_{d - 2} + \cdots + \xi_{1} + \xi_{0}\, \ge\, 1} \phi_{\left( \xi_{d - 2}, \cdots, \xi_{1}, \xi_{0} \right)} (z) A_{d - 2}^{\xi_{d - 2}} \cdots A_{1}^{\xi_{1}} A_{0}^{\xi_{0}}. \]
Then, substituting this expression for $\Phi$ in equation \eqref{conjugacy} yields 
\begin{eqnarray} 
\label{conjatwork}
& & z^{d} + \sum_{\xi_{d - 2} + \cdots + \xi_{1} + \xi_{0}\, \ge\, 1} \phi_{\left( \xi_{d - 2}, \cdots, \xi_{1}, \xi_{0} \right)} (z^{d}) A_{d - 2}^{\xi_{d - 2}} \cdots A_{1}^{\xi_{1}} A_{0}^{\xi_{0}} \nonumber \\ 
& - & \ \ \ \ \ \ \ \left( z + \sum_{\xi_{d - 2} + \cdots + \xi_{1} + \xi_{0}\, \ge\, 1} \phi_{\left( \xi_{d - 2}, \cdots, \xi_{1}, \xi_{0} \right)} (z) A_{d - 2}^{\xi_{d - 2}} \cdots A_{1}^{\xi_{1}} A_{0}^{\xi_{0}} \right)^{d} \nonumber \\ 
& - & A_{d - 2} \left( z + \sum_{\xi_{d - 2} + \cdots + \xi_{1} + \xi_{0}\, \ge\, 1} \phi_{\left( \xi_{d - 2}, \cdots, \xi_{1}, \xi_{0} \right)} (z) A_{d - 2}^{\xi_{d - 2}} \cdots A_{1}^{\xi_{1}} A_{0}^{\xi_{0}} \right)^{d - 2} \nonumber \\ 
& - & \cdots \nonumber \\ 
& - & A_{1} \ \ \ \left( z + \sum_{\xi_{d - 2} + \cdots + \xi_{1} + \xi_{0}\, \ge\, 1} \phi_{\left( \xi_{d - 2}, \cdots, \xi_{1}, \xi_{0} \right)} (z) A_{d - 2}^{\xi_{d - 2}} \cdots A_{1}^{\xi_{1}} A_{0}^{\xi_{0}} \right) - A_{0} \nonumber \\
= & & 0.
\end{eqnarray} 

\noindent 
Using the above equation, one can obtain the expression for all the functions $\phi_{\left( \xi_{d - 2}, \cdots, \xi_{1}, \xi_{0} \right)} (z)$ for the various $\xi_{r} \ge 0$ satisfying $\sum\limits_{r = 1}^{d - 2} \xi_{r} \ge 1$. However, since the theorems \eqref{theorem1} and \eqref{theorem2} only deal with the first and the second derivatives of the Lyapunov exponent, we are only interested in the cases when $\sum\limits_{r = 1}^{d - 2} \xi_{r} = 1$ and $\sum\limits_{r = 1}^{d - 2} \xi_{r} = 2$. The former case is achieved only when one and only one of the $\xi_{r} = 1$ and the remaining $\xi_{s} = 0$ while the latter case is achieved when either one and only one of the $\xi_{r} = 2$ and the remaining $\xi_{s} = 0$ or when two of the $\xi_{r} = 1$ and the remaining $\xi_{s} = 0$. 
\medskip 

\noindent 
For ease of writing and the readers' convenience, we make the following notations. 
\begin{eqnarray} 
\label{coefffn} 
\phi_{r} (z) & = & \phi_{\left( \xi_{d - 2}, \cdots, \xi_{1}, \xi_{0} \right)} (z)\ \text{when}\ \xi_{r} = 1\ \text{and}\ \xi_{s} = 0,\ \forall d - 2 \ge s \ge 0\ \text{with}\ s \ne r; \nonumber \\ 
\phi_{r^{2}} (z) & = & \phi_{\left( \xi_{d - 2}, \cdots, \xi_{1}, \xi_{0} \right)} (z)\ \text{when}\ \xi_{r} = 2\ \text{and}\ \xi_{s} = 0,\ \forall d - 2 \ge s \ge 0\ \text{with}\ s \ne r; \nonumber \\ 
\phi_{rs} (z) & = & \phi_{\left( \xi_{d - 2}, \cdots, \xi_{1}, \xi_{0} \right)} (z)\ \text{when}\ \xi_{r} = \xi_{s} = 1\ \text{with}\ r < s\ \text{and} \nonumber \\ 
& & \qquad \qquad \qquad \qquad \qquad \xi_{t} = 0,\ \forall d - 2 \ge t \ge 0\ \text{with}\ t \ne r\ \text{and}\ t \ne s. 
\end{eqnarray}

\noindent 
Thus, using these notations as in \eqref{coefffn} and making necessary computations as far as necessary in equation \eqref{conjatwork}, we obtain 
\begin{eqnarray}
\label{A_{r}}
\phi_{r} (z) & = & - z \sum_{\kappa_{1} \ge 1} \frac{1}{d^{\kappa_{1}}} \frac{1}{ z^{d^{\kappa_{1}} - r d^{\kappa_{1} - 1}}}\ ; \\ 
\label{A_{r}^{2}}
\phi_{r^{2}} (z) & = & -z \bigg[ \frac{d (d - 1)}{2} \sum_{\kappa_{3} \ge 1} \frac{1}{d^{\kappa_{3}}} \sum_{\kappa_{2} \ge 1} \sum_{\kappa_{1} \ge 1}^{\kappa_{2}} \frac{1}{d^{\kappa_{2} + 1}} \frac{1}{z^{d^{\kappa_{3}} - d^{\kappa_{3} - 1} + d^{\kappa_{3} - 1} (d^{\kappa_{1}} - r d^{\kappa_{1} - 1} + d^{\kappa_{2} - \kappa_{1} + 1} - r d^{\kappa_{2} - \kappa_{1}} - d + 1)}} \nonumber \\
& & \qquad \qquad - r \sum_{\kappa_{3} \ge 1} \frac{1}{d^{\kappa_{3}}} \sum_{\kappa_{1}\geq 1} \frac{1}{d^{\kappa_{1}}} \frac{1}{z^{d^{\kappa_{3}} - d^{\kappa_{3} - 1} + d^{\kappa_{3} - 1} (d^{\kappa_{1}} - r d^{\kappa_{1} - 1} - r + 1)}} \bigg]\ ; \\
 \label{A_{r}A_{s}}
\phi_{rs} (z) & = & -z \bigg[ d (d - 1) \sum_{\kappa_{3} \ge 1} \frac{1}{d^{\kappa_{3}}} \sum_{\kappa_{2} \ge 1} \sum_{\kappa_{1} \ge 1}^{\kappa_{2}} \frac{1}{d^{\kappa_{2} + 1}} \frac{1}{z^{d^{\kappa_{3}} - d^{\kappa_{3} - 1} + d^{\kappa_{3} - 1} (d^{\kappa_{1}} - r d^{\kappa_{1} - 1} + d^{\kappa_{2} - \kappa_{1} + 1} - s d^{\kappa_{2} - \kappa_{1}} - d + 1)}} \nonumber \\
& & \qquad \qquad - r \sum_{\kappa_{3} \ge 1} \frac{1}{d^{\kappa_{3}}} \sum_{\kappa_{1} \ge 1} \frac{1}{d^{\kappa_{1}}} \frac{1}{z^{d^{\kappa_{3}} - d^{\kappa_{3} - 1} + d^{\kappa_{3} - 1} (d^{\kappa_{1}} - s d^{\kappa_{1} - 1} - r + 1)}} \nonumber \\
& & \qquad \qquad - s \sum_{\kappa_{3} \ge 1} \frac{1}{d^{\kappa_{3}}} \sum_{\kappa_{1} \ge 1} \frac{1}{d^{\kappa_{1}}} \frac{1}{z^{d^{\kappa_{3}} - d^{\kappa_{3} - 1} + d^{\kappa_{3} - 1} (d^{\kappa_{1}} - r d^{\kappa_{1} - 1} - s + 1)}} \bigg]. 
\end{eqnarray}

 \section{Computation of Lyapunov Exponents}
\label{lecom} 

\noindent 
We know from equation \eqref{pratwork} that in order to calculate the Lyapunov exponent of $P$ with respect to the measure $\mu_{\vec{p}}$, it is sufficient for us to evaluate $\int_{\mathbb{S}^{1}} \log |\Phi_{P} (z)| d \mu_{\vec{p}}$. We now undertake the necessary computations, here. 
\begin{eqnarray} 
\label{congd}
- \int_{\mathbb{S}^{1}} \log | \Phi_{P} (z) | d \mu_{\vec{p}} & = & - \int_{\mathbb{S}^{1}} \sum_{r\, =\, 0}^{d - 2} {\rm Re} \left( A_{r} \overline{z} \phi_{r} \right) d \mu_{\vec{p}} \nonumber \\ 
& & - \int_{\mathbb{S}^{1}} \sum_{r\, =\, 0}^{d - 2} \left[ {\rm Re} \left( A_{r}^{2} \overline{z} \phi_{r^{2}} \right) - \frac{1}{2} \left\{{\rm Re} \left( A_{r} \overline{z} \phi_{r} \right) \right\}^{2} + \frac{1}{2} \left\{{\rm Im} \left( A_{r} \overline{z} \phi_{r} \right) \right\}^{2} \right] d \mu_{\vec{p}} \nonumber \\  
& & - \int_{\mathbb{S}^{1}} \sum_{r\, =\, 0}^{d - 3}\ \sum_{r\, <\, s\, =\, 1}^{d - 2}\Big[ {\rm Re} \left( A_{r} A_{s} \overline{z} \phi_{rs} \right) - \left\{{\rm Re} \left( A_{r}\overline{z} \phi_{r} \right) \times {\rm Re} \left( A_{s} \overline{z} \phi_{s} \right) \right\} \nonumber \\ 
& & \qquad \qquad \qquad \qquad \qquad \qquad + \left\{{\rm Im} \left( A_{r} \overline{z} \phi_{r} \right) \times {\rm Im} \left( A_{s} \overline{z} \phi_{s} \right) \right\} \Big] d \mu_{\vec{p}} \nonumber \\ 
& & +\ O \left( \text{terms where}\ \sum_{r\, =\, 0}^{d - 2} \xi_{r} \ge 3 \right). 
\end{eqnarray} 

\noindent 
We first observe that each of the integrals in the right hand side of equation \eqref{congd} evaluates to $0$ when the measure of integration is the equidistributed Lyubich's measure, as defined in section \eqref{weight}. We now evaluate the integrals with respect to the weighted Lyubich's measure, also defined in section \eqref{weight}, in particular when one of the $p_{j} \uparrow 1$ for some $d \ge j \ge 1$. We separate the cases of the polynomial, $P_{\mathbb{R}}$, as written in equation \eqref{polyreal} with only real coefficients, \textit{i.e.}, $A_{r} = \alpha_{r}$ and $\beta_{r} \equiv 0\ \forall d - 2 \ge r \ge 0$ and $P_{\mathbb{C}}$, as written in equation \eqref{polycomp} with complex coefficients, \textit{i.e.}, $A_{r} = \alpha_{r} + i \beta_{r},\ \forall d - 2 \ge r \ge 0$. 

\subsection*{Computations for $P_{\mathbb{R}}$}

\noindent 
We urge the reader to observe that when the coefficients of the polynomial are all real, we take the terms containing ${\rm Re} A_{r} = \alpha_{r}$ out of the integral, as a multiplicative factor. Further, in this case, we have ${\rm Im} A_{r} = 0$. Thus, 
\begin{eqnarray*} 
\int {\rm Re} \left( \overline{z} \phi_{r} \right) d \mu_{\vec{p}} & \to & - \frac{1}{d - 1}\ ; \\ 
\\ 
\int {\rm Re} \left( \overline{z} \phi_{r^{2}} \right) d \mu_{\vec{p}} & \to & - \frac{d - 2r}{2 (d - 1)^{2}}\ ; \\ 
\\ 
\frac{1}{2} \int \left[ {\rm Re} \left( \overline{z} \phi_{r} \right) \right]^{2} d \mu_{\vec{p}} & \to & \frac{1}{2} \frac{1}{(d - 1)^{2}}\ ; \\ 
\\ 
\int {\rm Re} \left( \overline{z} \phi_{rs} \right) d \mu_{\vec{p}} & \to & - \frac{d - r - s}{(d - 1)^{2}}\ ; \\ 
\\ 
\int \left[ {\rm Re} \left( \overline{z} \phi_{r} \right) \times {\rm Re} \left( \overline{z} \phi_{s} \right) \right] d \mu_{\vec{p}} & \to & \frac{1}{(d - 1)^{2}}. 
\end{eqnarray*} 

\subsection*{Computations for $P_{\mathbb{C}}$} 

\noindent 
In this case, $A_{r} = \alpha_{r} + i \beta_{r}$ for all $d - 2 \ge r \ge 0$. Then, the computations yield 
\begin{eqnarray*} 
\int {\rm Re} \left( A_{r} \overline{z} \phi_{r} \right) d \mu_{\vec{p}} & \to & - \frac{1}{d - 1} \alpha_{r}\ ; \\ 
\\ 
\int {\rm Re} \left( A_{r}^{2} \overline{z} \phi_{r^{2}} \right) d \mu_{\vec{p}} & \to & - \frac{d - 2r}{2 (d - 1)^{2}}\left( \alpha_{r}^{2} - \beta_{r}^{2} \right)\ ; \\ 
\\ 
\frac{1}{2} \int \left[ {\rm Re} \left( A_{r} \overline{z} \phi_{r} \right) \right]^{2} d \mu_{\vec{p}} & \to & \frac{1}{2} \frac{1}{(d - 1)^{2}} \alpha_{r}^{2}\ ; \\ 
\\ 
\frac{1}{2} \int \left[ {\rm Im} \left( A_{r} \overline{z} \phi_{r} \right) \right]^{2} d \mu_{\vec{p}} & \to & \frac{1}{2} \frac{1}{(d - 1)^{2}} \beta_{r}^{2}\ ; \\ 
\\ 
\int {\rm Re} \left( A_{r} A_{s} \overline{z} \phi_{rs} \right) d \mu_{\vec{p}} & \to & - \frac{d - r - s}{(d - 1)^{2}} \left( \alpha_{r} \alpha_{s} - \beta_{r} \beta_{s} \right)\ ; \\
\\ 
\int \left[ {\rm Re} \left( A_{r} \overline{z} \phi_{r} \right) \times {\rm Re} \left( A_{s} \overline{z} \phi_{s} \right) \right] d \mu_{\vec{p}} & \to & \frac{1}{(d - 1)^{2}} \left( \alpha_{r} \alpha_{s} \right)\ ; \\ 
\\ 
\int \left[ {\rm Im} \left( A_{r} \overline{z} \phi_{r} \right) \times {\rm Im} \left( A_{s} \overline{z} \phi_{s} \right) \right] d \mu_{\vec{p}} & \to & \frac{1}{(d - 1)^{2}} \left( \beta_{r} \beta_{s} \right). 
\end{eqnarray*} 

\noindent 
Thus, from the computations of the Lyapunov exponent, we have, as $p_{j} \uparrow 1$, for some $d \ge j \ge 1$, that  
\begin{eqnarray} 
\label{lyapreal} 
\Lambda_{\mu_{\vec{p}}} (P_{\mathbb{R}}) & \to & - \log d\ +\ \sum_{r\, =\, 0}^{d - 2} \frac{1}{d - 1} \alpha_{r}\ +\ \sum_{r\, =\, 0}^{d - 2} \frac{d - 2r + 1}{2 (d - 1)^{2}} \alpha_{r}^{2} \nonumber \\ 
& & +\ \sum_{r\, =\, 0}^{d - 3}\ \sum_{r\, <\, s\, =\, 1}^{d - 2} \frac{d - r - s + 1}{(d - 1)^{2}} \alpha_{r} \alpha_{s}. \\
\label{lyapcomplex}
\Lambda_{\mu_{\vec{p}}} (P_{\mathbb{C}}) & \to & - \log d\ +\ \sum_{r\, =\, 0}^{d - 2} \frac{1}{d - 1} \alpha_{r}\ +\ \sum_{r\, =\, 0}^{d - 2} \frac{d - 2r + 1}{2 (d - 1)^{2}} \alpha_{r}^{2}\ -\ \sum_{r\, =\, 0}^{d - 2} \frac{d - 2r + 1}{2 (d - 1)^{2}} \beta_{r}^{2} \nonumber \\
& & +\ \sum_{r\, =\, 0}^{d - 3}\ \sum_{r\, <\, s\, =\, 1}^{d - 2} \frac{d - r - s + 1}{(d - 1)^{2}} \alpha_{r} \alpha_{s}\ -\ \sum_{r\, =\, 0}^{d - 3}\ \sum_{r\, <\, s\, =\, 1}^{d - 2} \frac{d - r - s + 1}{(d - 1)^{2}} \beta_{r} \beta_{s}.
\end{eqnarray} 

\section{Proofs of the main theorems and corollaries} 
\label{proof}

\noindent 
In this concluding section, we write the proof of the main theorems \eqref{theorem1} and \eqref{theorem2}, by appealing to the computations that we have done in section \eqref{lecom}. 
\medskip 

\noindent
\begin{proof} (of theorem \eqref{theorem1}) 
We obtain the following expressions, by directly differentiating equation \eqref{lyapcomplex}, with respect to the corresponding variable. 
\begin{eqnarray*} 
\frac{\partial^{2} \Lambda_{\mu_{\vec{p}}}}{\partial \alpha_{r}^{2}}\ =\ \frac{d - 2r + 1}{2 (d - 1)^{2}}, & & \frac{\partial^{2} \Lambda_{\mu_{\vec{p}}}}{\partial \beta_{r}^{2}}\ =\ - \frac{d - 2r + 1}{2 (d - 1)^{2}}, \\ 
\\ 
\frac{\partial^{2} \Lambda_{\mu_{\vec{p}}}}{\partial \alpha_{r}{\partial \alpha_{s}}}\ =\ \frac{d - r - s + 1}{(d - 1)^{2}}, & & \frac{\partial^{2} \Lambda_{\mu_{\vec{p}}}}{\partial \beta_{r}{\partial \beta_{s}}}\ =\ -\frac{d - r - s + 1}{(d - 1)^{2}}\ \ \text{for}\ d - 2 \ge s > r \ge 0. 
\end{eqnarray*} 
Thus, we have 
\[ \frac{\partial^{2} \Lambda_{\mu_{\vec{p}}}}{\partial \alpha_{r}{\partial \alpha_{s}}}\ \ =\ \ - \frac{\partial^{2} \Lambda_{\mu_{\vec{p}}}}{\partial \beta_{r}{\partial \beta_{s}}},\ \ \ \forall d - 2 \ge s > r \ge 0. \]
\end{proof} 
\medskip 

\noindent 
\begin{proof} (of theorem \eqref{theorem2}) 
We obtain the following expressions, by directly differentiating equation \eqref{lyapcomplex} and \eqref{lyapreal}, with respect to the corresponding variable. 
\begin{eqnarray*} 
\frac{\partial \Lambda_{\mu} (P_{\mathbb{R}})}{\partial \alpha_{r}} & = & \frac{\partial \Lambda_{\mu} (P_{\mathbb{C}})}{\partial \alpha_{r}} = \frac{1}{d - 1} + \frac{d - 2r + 1}{2(d - 1)^{2}} \alpha_{r} + \sum_{r\, \ne\, s\, =\, 0}^{d - 2} \frac{d - r - s}{(d - 1)^{2}} \alpha_{s} ; \\ 
\frac{\partial^{2} \Lambda_{\mu} (P_{\mathbb{R}})}{\partial \alpha_{r}^{2}} & = & \frac{\partial^{2} \Lambda_{\mu} (P_{\mathbb{C}})}{\partial \alpha_{r}^{2}} = \frac{d - 2r + 1}{2 (d - 1)^{2}}. 
\end{eqnarray*} 
\end{proof} 
\medskip 

\noindent 
The observations in \cite{SridharanArt} by the authors, are then simple corollaries to specific cases when $d = 2$ and $d = 3$. We complete this paper with these following corollaries. 
\medskip 

\noindent 
\begin{corollary} 
For a monic centred hyperbolic quadratic polynomial 
\[ R_{2} (z)\ \ =\ \ z^{2} + \alpha + i \beta\ \ \ \ \text{with}\ \ \alpha^{2} + \beta^{2} < 1, \] 
we have 
\[ \Lambda_{\mu_{\vec{p}}} (R_{2})\ \ \to\ \ - \log 2\ +\ \alpha\ +\ \frac{3}{2} \left( \alpha^{2} - \beta^{2} \right),\ \ \text{as}\ p_{1} \uparrow 1\ \text{or as}\ p_{2} \uparrow 1. \] 
\end{corollary} 
\medskip 

\noindent 
\begin{corollary} 
For a monic centred hyperbolic cubic polynomial 
\[ R_{3} (z)\ \ =\ \ z^{3} + \left( \alpha_{1} + i \beta_{1} \right) z + \left( \alpha_{0} + i \beta_{0} \right)\ \ \ \ \text{with}\ \ \alpha_{r}^{2} + \beta_{r}^{2} < 1\ \ \text{for}\ r = 1, 0, \] 
we have 
\begin{eqnarray*} 
\Lambda_{\mu_{\vec{p}}} (R_{3}) & \to & - \log 3\ +\ \frac{1}{2} \left( \alpha_{1} + \alpha_{0} \right)\ +\ \frac{1}{2} \left( \alpha_{0}^{2} - \beta_{0}^{2} \right)\ +\ \frac{1}{4} \left( \alpha_{1}^{2} - \beta_{1}^{2} \right)\ +\ \frac{3}{4} \left( \alpha_{1} \alpha_{0} - \beta_{1} \beta_{0} \right) \\ 
& & \hspace{+8cm} \text{as}\ p_{1} \uparrow 1,\ p_{2} \uparrow 1\ \text{or as}\ p_{3} \uparrow 1. 
\end{eqnarray*} 
\end{corollary} 

\bigskip 

\noindent 
{\sc Sridharan, Shrihari} \\
Indian Institute of Science Education and Research Thiruvananthapuram (IISER-TVM) \\
Maruthamala P.O., Vithura, Thiruvananthapuram, INDIA. PIN 695 551. \\ 
{\tt shrihari@iisertvm.ac.in}  
\bigskip 

\noindent 
{\sc Tiwari, Atma Ram} \\ 
Indian Institute of Science Education and Research Thiruvananthapuram (IISER-TVM) \\
Maruthamala P.O., Vithura, Thiruvananthapuram, INDIA. PIN 695 551. \\ 
{\tt artiwari15@iisertvm.ac.in}

\bigskip

\noindent
\begin{thebibliography}{99}

\bibitem{bcnv:18} 
{\sc Barabanov, E.}, {\sc Czornik, A.}, {\sc Niezabitowski, M.} and {\sc Vaidzelevich, A.}, ``Influence of parametric perturbations on Lyapunov exponents of discrete linear time-varying systems", \emph{Systems Control Lett.}, {\bf 122}, (2018), 54 - 59. 

\bibitem{Beardon}
{\sc Beardon, A. F.}, \emph{Iteration of rational functions: Complex analytic dynamical systems}, Graduate Texts in Mathematics, {\bf 132}, Springer-Verlag, New York, (1991). 

\bibitem{Benettin}
{\sc Benettin, G.}, {\sc Galgani, L.} and {\sc Strelcyn, J.-M.}, ``Kolmogorov entropy and numerical experiments" (1976), \emph{Physical Review A}, 14 (1976), 2338 - 2345.

\bibitem{cg:93} 
{\sc Carleson, L.} and {\sc Gamelin, T.W.}, \emph{Complex Dynamics}, Universitext: Tracts in Mathematics, Springer-Verlag, New York, 1993. 

\bibitem{tc:19} 
{\sc Catalan, T.}, ``A link between topological entropy and Lyapunov exponents", \emph{Ergod. Th. Dynam. Sys.}, {\bf 39}, (2019), 620 - 637. 

\bibitem{Coelho}
{\sc Coelho, Z.} and {\sc Parry, W.}, ``Central limit asymptotics for shifts of finite type", \emph{Israel J. Math.}, {\bf 69}, (1990), 235 - 249.

\bibitem{gr:94} 
{\sc Gozzi, E.} and {\sc Reuter, M.}, ``Lyapunov exponents, path-integrals and forms", \emph{Chaos, solitons and fractals}, {\bf 4}, (1994), 1117 - 1139.  

\bibitem{Lyubich} 
{\sc Lyubich, M. Yu.}, ``The dynamics of rational transforms: the topological picture", \emph{Russian Math. Surveys}, {\bf 41}, (1986), 43 - 117.

\bibitem{Pesin}
{\sc Pesin, Ya. B.}, ``Characteristic Lyapunov exponents and smooth ergodic theory." \emph{Russian Math. Surveys}, {\bf 32}, (1977), 55 - 114.

\bibitem{Ruelle}
{\sc Ruelle, D.}, \emph{Thermodynamic formalism}, Encyclopedia Mathematics and its Applications, Reading: Addison-Wesley, (1978).

\bibitem{Sridharan} 
{\sc Sridharan, S.}, ``Non-vanishing derivatives of Lyapunov exponents and the pressure function", \emph{Dyn. Syst.}, {\bf 21}, (2006), 491 - 500.

\bibitem{SridharanArt} 
{\sc Sridharan, S.} and {\sc Tiwari, A.,R.}, ``The dependence of Lyapunov exponents of polynomials on their coefficients", \emph{J. Comput. Dyn.}, {\bf 6}, (2019), 95 - 109.

\bibitem{Steinmetz}
{\sc Steinmetz, N.}, \emph{Rational iteration: complex analytic dynamical systems}, De Gruyter studies in Mathematics, {\bf 16}, Walter de Gruyter and Co., Berlin, (1993).

\bibitem{Walter}
{\sc Walters, P.}, \emph{An introduction to ergodic theory}, Graduate texts in Mathematics, {\bf 79}, Springer-Verlag, New York, (1982).

\bibitem{Zinsmeister}
{\sc Zinsmeister, M.}, \emph{Formalisme thermodynamique et syst\'{e}mes dynamiques holomorphes}", Panoramas et Synth\'{e}ses, {\bf 4}, (1996).

\end{thebibliography}
\end{document}